# Energy Storage Sharing Strategy in Distribution Networks Using Bi-level Optimization Approach

Huimiao Chen, *Student Member, IEEE*, Yang Yu, Zechun Hu, *Member, IEEE*, Haocheng Luo, Chin-Woo Tan, Ram Rajagopal, *Member, IEEE*

*Abstract*—In this paper, we address the energy storage management problem in distribution networks from the perspective of an independent energy storage manager (IESM) who aims to realize optimal energy storage sharing with multi-objective optimization, i.e., optimizing the system peak loads and the electricity purchase costs of the distribution company (DisCo) and its customers. To achieve the goal of the IESM, an energy storage sharing strategy is therefore proposed, which allows DisCo and customers to control the assigned energy storage. The strategy is updated day by day according to the system information change. The problem is formulated as a bi-level mathematical model where the upper level model (ULM) seeks for optimal division of energy storage among Disco and customers, and the lower level models (LLMs) represent the minimizations of the electricity purchase costs of DisCo and customers. Further, in order to enhance the computation efficiency, we transform the bi-level model into a single-level mathematical program with equilibrium constraints (MPEC) model and linearize it. Finally, we validate the effectiveness of the strategy and complement our analysis through case studies.

*Index Terms*—Energy storage sharing, scheduling, demand response, MPEC, smart grids.

## I. INTRODUCTION

SMART grid development calls for effective solutions, such as flexible loads and energy storage systems (ESSs), to facilitate demand response and meet the energy and environmental challenges. Compared with flexible loads, such as electric vehicles (EVs) and washing machines, energy storage provides a better quality of flexibility in power systems since it can store energy when the electricity supply is surplus and release energy during the high demand period to alleviate network congestion [1], [2]. Additionally, as the ESS technology advances [3] and its costs drop [2], a number of energy storage projects in power systems have been started up or been placed on the agenda all over the world. For example, in China, a peaking power plant based on the technology of chemical energy storage will be built in the northeast city of Dalian, aiming to reach a record capacity of 200 MW upon completion in 2018 [4].

For customers facing time-varying prices, energy storage provides a means to reduce energy costs through arbitrage, i.e., applying the principle of "buy low, sell high". References [5]-[7] investigate the optimal energy storage operation policy for a single customer to maximize the its benefit. For energy storage owned by network operators, reference [8] presents an asymptotically optimal control policy for large-scale energy storage and authors of [9] propose a dynamic program-based method to estimate the capacity value of energy storage.

Further, under the energy storage sharing scenario, a large-scale ESS can not only mitigate the economic deterrents of each customer caused by expensive purchase and maintenance costs of small-scale energy storage devices, but also provide a more manageable platform for demand response. Along with these advantages, sharing makes the challenges of energy storage management and control increase. In the existing literature, authors of [10] use a Markovian model to address the cost saving trade-off problem of sharing ESS among a group of customers in a community. Reference [11] discusses the energy storage managing method in distribution network based on evenly dividing energy storage between customers and system operator, but does not optimize the division of energy storage. In [2], energy storage is also shared between customers and system operator, and customers seek for lower wholesale energy costs, whereas system operator aims to minimize network investment costs. However, in [2], only five defined energy storage dispatch scenarios are used as the candidates for the sharing, which leads to suboptimal results.

In this paper, for the purpose of achieving the optimal storage management among the distribution company (DisCo) and customers, we put forward an energy storage sharing strategy from the perspective of an independent energy storage manager (IESM), who is responsible for the day-ahead energy storage division. The IESM is able to obtain the system information, such as locational marginal prices (LMPs) and customer load forecasting data. Then, a bi-level optimization model is built to solve the problem, where the upper level model (ULM) deals with the storage division problem for the optimization of the social welfare, i.e., the system peak loads and the electricity purchase costs of DisCo and customers, whereas the lower level models (LLM) describes Disco's and customers' behaviors of using a given capacity of energy storage for the minimum energy costs. To effectively solve the model, the complementarity theory is utilized to tackle the bi-level model and convert it into a single-level mathematical program with equilibrium constraints (MPEC) model.

This work was supported in part by the National Natural Science Foundation of China under Grant 51477082.

H. Chen, Z. Hu and H. Luo are with the Department of Electrical Engineering, Tsinghua University, Beijing, 100084, P. R. China (email: chenhm15@mails.tsinghua.edu.cn).

Y. Yu, C. Tan and R. Rajagopal are with the Department of Civil and Environmental Engineering, Stanford University, Stanford, CA, 94305, USA.

Additionally, the complementarity constraints in the MPEC model are linearized, and then the model becomes a mixed integer linear programming (MILP) problem, which can be easily solved via mature optimization software. Finally, the effectiveness of the proposed strategy is verified and analyzed in case studies.

## II. MATHEMATICAL FORMULATION OF ENERGY STORAGE SHARING STRATEGY

### A. Strategy Overview

In the scenario of the proposed energy storage sharing strategy, the IESM divides the energy storage with given capacity among a DisCo and its customers at day-ahead period. For preparation, all the customer need to report their predicted electricity usage for the next day. Besides, LMPs, i.e., the prices cleared by the independent system operator (ISO) in the day-ahead electricity market, and time-of-use (TOU) prices are adopted as the electricity purchase prices for the DisCo and the customers, respectively. Thus, the energy storage sharing problem is solved after the day-head electricity market clearing. The final energy storage sharing decisions will be carried out at the beginning of the next day and then the new cycle starts (see Fig. 1).

### B. Bi-level Model Formulation

In this section, a bi-level model is formulated to achieve the optimal energy storage sharing strategy, which includes the energy storage division of the IESM in ULM and the energy storage operation control of the DisCo and customers in LLMs. In the bi-level model, LLMs are actually the constraints of the ULM.

*1) ULM: Energy Storage Division of IESM*

In the ULM, we suppose that the IESM aims to i) shave the system peak load $P^{peak}$ so as to defer the network investment; ii) optimize the energy purchase costs of both DisCo, i.e., $C^d$, and customers, i.e., $C^c$. Note that the objective can be altered according to the needs in practical operation. The ULM is modeled as (1-a)-(1-e).

$$\min \lambda_1 \cdot P^{peak} + \lambda_2 \cdot C^d + \lambda_3 \cdot C^c \quad (1\text{-a})$$

where

$$C^d = \sum_{t=1}^{T} \pi_t^{LMP} \cdot \left( P_t^{ori,d} + \sum_{n=1}^{N} P_{n,t}^{ch,c} - \sum_{n=1}^{N} P_{n,t}^{dis,c} + P_t^{ch,d} - P_t^{dis,d} \right) \cdot \Delta t \quad (1\text{-b})$$

$$C^c = \sum_{t=1}^{T} \pi_t^{TOU} \cdot \left( P_t^{ori,d} + \sum_{n=1}^{N} P_{n,t}^{ch,c} - \sum_{n=1}^{N} P_{n,t}^{dis,c} + P_t^{ch,d} - P_t^{dis,d} \right) \cdot \Delta t \quad (1\text{-c})$$

subject to:

$$S^{ES,d} + \sum_{n=1}^{N} S_n^{ES,c} \leq S^{ES,total} \quad (1\text{-d})$$

$$P^{peak} \geq P_t^{ori,d} + \sum_{n=1}^{N} P_{n,t}^{ch,c} - \sum_{n=1}^{N} P_{n,t}^{dis,c} + P_t^{ch,d} - P_t^{dis,d}, \forall t. \quad (1\text{-e})$$

In (1-a), $\lambda_1$, $\lambda_2$ and $\lambda_3$ are weight factors. And in (1-b) and (1-c), $\pi_t^{LMP}$ and $\pi_t^{TOU}$ are respectively LMP and TOU price at time slot $t$ (suppose one day is divided equally into $T$ discrete time slots), and $P_t^{ori,d}$, $P_{n,t}^{ch,c}$, $P_{n,t}^{dis,c}$, $P_t^{ch,d}$ and $P_t^{dis,d}$ are the original system load at time slot $t$, the charging power and discharging power of customer $n$ at time slot $t$, and charging power and discharging power of DisCo at time slot $t$,

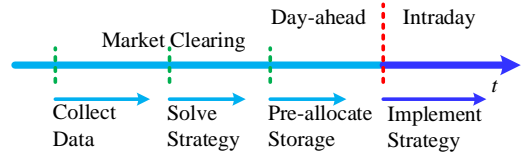

Fig. 1 The schematic diagram of strategy implementation steps in a daily cycle.

respectively. $N$ is the total number of customers and $\Delta t$ is the duration of a time slot.

In the constraints of ULM, (1-d) is the capacity constraint of energy storage, where $S^{ES,d}$, $S_n^{ES,c}$ and $S^{ES,total}$ denote the shared energy storage capacity for DisCo and customer $n$, and the total capacity, respectively; (1-e) are linearized expressions to describe $P^{peak}$, because the optimal solution must meet equation (2) to achieve the lowest peak load.

$$P^{peak} = \max_{\forall t} \left( P_t^{ori,d} + \sum_{n=1}^{N} P_{n,t}^{ch,c} - \sum_{n=1}^{N} P_{n,t}^{dis,c} + P_t^{ch,d} - P_t^{dis,d} \right) \quad (2)$$

In the ULM, $P^{peak}$, $S^{ES,d}$ and all the $S_n^{ES,c}$ compose the decision variables.

*2) LLM-C: Energy Storage Operation Control of Customers*

In the LLM of customer $n$ (LLM-C), the energy storage capacity $S_n^{ES,c}$ becomes the input parameter and the customer tries to reduce his or her electricity bill. The model is formulated as (3-a)-(3-f).

$$\min \sum_{t=1}^{T} \pi_t^{TOU} \cdot \left( P_{n,t}^{ch,c} - P_{n,t}^{dis,c} \right) \cdot \Delta t + \alpha \cdot \left( P_n^{peak,c} - P_n^{valley,c} \right) \quad (3\text{-a})$$

subject to:

$$S_n^{ES,c} \cdot SoC^{lower} \leq S_n^{ES,c} \cdot SoC_n^{ini,c} + \sum_{\tau=1}^{t} P_{n,\tau}^{ch,c} \cdot \Delta t \cdot \eta^{ch}$$
$$- \sum_{\tau=1}^{t} P_{n,\tau}^{dis,c} \cdot \Delta t / \eta^{dis} \leq S_n^{ES,c} \cdot SoC^{upper}, \forall t \quad (3\text{-b})$$

$$\sum_{t=1}^{T} P_{n,t}^{ch,c} \cdot \Delta t \cdot \eta^{ch} - \sum_{t=1}^{T} P_{n,t}^{dis,c} \cdot \Delta t / \eta^{dis} = 0 \quad (3\text{-c})$$

$$P_n^{peak,c} \geq P_{n,t}^{ori,c} + P_{n,t}^{ch,c} - P_{n,t}^{dis,c}, \forall t \quad (3\text{-d})$$

$$P_n^{valley,c} \leq P_{n,t}^{ori,c} + P_{n,t}^{ch,c} - P_{n,t}^{dis,c}, \forall t \quad (3\text{-e})$$

$$0 \leq P_{n,t}^{dis,c}, P_{n,t}^{ch,c} \leq k \cdot S_n^{ES,c}, \forall t. \quad (3\text{-f})$$

In the objective of LLM-C, the first term is the increment of electricity purchase costs of customer $n$; the second term is a penalty term of peak-valley difference, i.e., $P_n^{peak,c} - P_n^{valley,c}$, with a small weight factor $\alpha$ because the energy storage operation problem usually has multiple optimal solutions under the TOU prices which only change one time several hours (the customer can charge or discharge the energy storage casually at the periods of same price without economic loss), and a small penalty term can make the solution unique and improve the performance with only a little sacrifice of the customer benefits. In practice, the penalty term can be regard as the extra payment to participate the energy storage sharing plan and can be altered according to different needs. For constraints, (3-b) restrict the states of charge (SoCs) of energy within the lower and upper bounds, i.e., $SoC^{lower}$ and $SoC^{upper}$, and $SoC_n^{ini,c}$ is the initial SoC of the energy storage controlled by customer $n$, and $\eta^{ch} / \eta^{dis}$ is charging/discharging efficiency of energy storage; (3-c) describes that energy charged to energy storage is equal to the one consumed since it is usually expected that the initial and final SoCs are same [12]; (3-d) and (3-e) are similar to (1-e) in the ULM, and $P_{n,t}^{\text{ori,c}}$ is the original load of customer $n$ at time slot $t$; (3-f) ensure that the charging power and discharging power are within the available range and the upper bound, i.e., $k \cdot S_n^{\text{ES,c}}$ where $k$ is a fixed coefficient, is proportional to the capacity of the allocated energy storage.

The decision variables of LLM-C are composed of $P_{n,t}^{\text{dis,c}}$, $P_{n,t}^{\text{ch,c}}$, $P_n^{\text{peak,c}}$ and $P_n^{\text{valley,c}}$.

*3) LLM-D: Energy Storage Operation Control of DisCo*

The LLM of DisCo (LLM-D) is quite similar to LLM-C, but the penalty term is not added here because different from TOU prices, the LMP $\pi_t^{\text{LMP}}$ varies frequently and the unique optimal energy storage operation solution generally exists. Note that a penalty term is allowed in practice if necessary. The model formulation is as follows:

$$\min \sum_{t=1}^{T} \pi_t^{\text{LMP}} \cdot \left( P_t^{\text{ch,d}} - P_t^{\text{dis,d}} \right) \cdot \Delta t \quad (4\text{-a})$$

subject to:

$$S^{\text{ES,d}} \cdot SoC^{\text{lower}} \leq S^{\text{ES,d}} \cdot SoC^{\text{ini,d}} + \sum_{\tau=1}^{t} P_\tau^{\text{ch,d}} \cdot \Delta t \cdot \eta^{\text{ch}}$$

$$-\sum_{\tau=1}^{t} P_\tau^{\text{dis,d}} \cdot \Delta t / \eta^{\text{dis}} \leq S^{\text{ES,d}} \cdot SoC^{\text{upper}}, \forall t \quad (4\text{-b})$$

$$\sum_{t=1}^{T} P_t^{\text{ch,d}} \cdot \Delta t \cdot \eta^{\text{ch}} - \sum_{t=1}^{T} P_t^{\text{dis,d}} \cdot \Delta t / \eta^{\text{dis}} = 0 \quad (4\text{-c})$$

$$0 \leq P_t^{\text{dis,d}}, P_t^{\text{ch,d}} \leq k \cdot S^{\text{ES,d}}, \forall t. \quad (4\text{-d})$$

In the above model, the objective is to minimize the increment of energy costs of DisCo, and constraints (4-b)-(4-d) are similar to constraints (3-b), (3-c) and (3-f), respectively.

The decision variables of LLM-D include $P_t^{\text{dis,d}}$ and $P_t^{\text{ch,d}}$.

## III. SOLUTION METHOD

### A. MPEC Model Formulation

Here, the complementarity theory used to tackle the proposed bi-level model. Some researchers have used this approach to solve the bi-level model-based bidding problem in electricity market [13], [14].

Since DisCo and customers take $S^{\text{ES,d}}$ and $S_n^{\text{ES,c}}$ as the parameters, all the LLMs are linear. Thus, they can be represented by their equivalent sets of Karush-Kuhn-Tucker (KKT) conditions and included as equality and inequality constraints of the ULM, then the bi-level model is transformed into a single-level MPEC model (see Fig. 2).

We use $f_n^c(x_n^c)$ and $f^d(x^d)$ to respectively denote the objective functions of customer $n$ and DisCo, where $x_n^c$ and $x^d$ are vectors composed of decision variables. Let $SoC_{n,t}^{\text{cur,c}}$ and $SoC_t^{\text{cur,d}}$ denote the current SoC of the energy storage respectively controlled by customer $n$ and DisCo at time slot $t$. Then, let:

i) $g_{n,1,t}^c(x_n^c) = SoC_{n,t}^{\text{cur,c}} - SoC^{\text{lower}}$, $g_{1,t}^d(x^d) = SoC_t^{\text{cur,d}} - SoC^{\text{lower}}$;

ii) $g_{n,2,t}^c(x_n^c) = SoC^{\text{upper}} - SoC_{n,t}^{\text{cur,c}}$, $g_{2,t}^d(x^d) = SoC^{\text{upper}} - SoC_t^{\text{cur,d}}$;

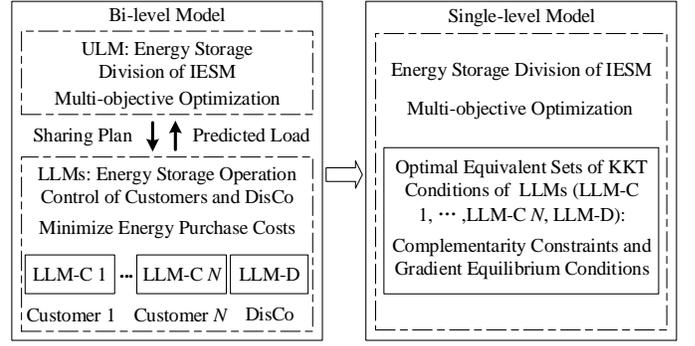

Fig. 2 Transformation from bi-level model to single-level model.

iii) $g_{n,3,t}^c(x_n^c) = P_{n,t}^{\text{dis,c}}$, $g_{3,t}^d(x^d) = P_t^{\text{dis,d}}$;

iv) $g_{n,4,t}^c(x_n^c) = P_{n,t}^{\text{ch,c}}$, $g_{4,t}^d(x^d) = P_t^{\text{ch,d}}$;

v) $g_{n,5,t}^c(x_n^c) = k \cdot S_n^{\text{ES,c}} - P_{n,t}^{\text{dis,c}}$, $g_{5,t}^d(x^d) = k \cdot S^{\text{ES,d}} - P_t^{\text{dis,d}}$;

vi) $g_{n,6,t}^c(x_n^c) = k \cdot S_n^{\text{ES,c}} - P_{n,t}^{\text{ch,c}}$, $g_{6,t}^d(x^d) = k \cdot S^{\text{ES,d}} - P_t^{\text{ch,d}}$;

vii) $g_{n,7,t}^c(x_n^c) = P_n^{\text{peak,c}} - P_{n,t}^{\text{ori,c}} - P_{n,t}^{\text{ch,c}} + P_{n,t}^{\text{dis,c}}$;

viii) $g_{n,8,t}^c(x_n^c) = P_{n,t}^{\text{ori,c}} + P_{n,t}^{\text{ch,c}} - P_{n,t}^{\text{dis,c}} - P_n^{\text{valley,c}}$;

ix) $h_n^c(x_n^c) = \sum_{t=1}^{T} P_{n,t}^{\text{ch,c}} \cdot \Delta t \cdot \eta^{\text{ch}} - \sum_{t=1}^{T} P_{n,t}^{\text{dis,c}} \cdot \Delta t / \eta^{\text{dis}}$,

$h^d(x^d) = \sum_{t=1}^{T} P_t^{\text{ch,d}} \cdot \Delta t \cdot \eta^{\text{ch}} - \sum_{t=1}^{T} P_t^{\text{dis,d}} \cdot \Delta t / \eta^{\text{dis}}$.

Then, the equivalent set of KKT conditions of LLM-C is as follows:

$$(3\text{-b})\text{-}(3\text{-f}) \quad (5\text{-a})$$

$$\nabla_{x_n^c} L_n^c(x_n^c, \omega_n^c, v_n^c) = \mathbf{0} \quad (5\text{-b})$$

$$\omega_{n,i,t}^c \cdot g_{n,i,t}^c(x_n^c) = 0, i = 1, 2, \cdots, 8, \forall t \quad (5\text{-c})$$

$$\omega_{n,i,t}^c \geq 0, i = 1, 2, \cdots, 8, \forall t \quad (5\text{-d})$$

where $L_n^c(x_n^c, \omega_n^c, v_n^c) = f_n^c(x_n^c) - \sum_{i=1}^{7}\sum_{t=1}^{T} \omega_{n,i,t}^c \cdot g_{n,i,t}^c(x_n^c) - v_n^c \cdot h_n^c(x_n^c)$, $\omega_n^c$ is the vector composed of $\omega_{n,i,t}^c$. $\omega_{n,i,t}^c$ and $v_n^c$ are dual variables. Similarly, the equivalent set of KKT conditions of LLM-D is as follows:

$$(4\text{-b})\text{-}(4\text{-d}) \quad (6\text{-a})$$

$$\nabla_{x^d} L^d(x^d, \omega^d, v^d) = \mathbf{0} \quad (6\text{-b})$$

$$\omega_{i,t}^d \cdot g_{i,t}^d(x^d) = 0, i = 1, 2, \cdots, 6, \forall t \quad (6\text{-c})$$

$$\omega_{i,t}^d \geq 0, i = 1, 2, \cdots, 6, \forall t \quad (6\text{-d})$$

where $L^d(x^d, \omega^d, v^d) = f^d(x^d) - \sum_{i=1}^{5}\sum_{t=1}^{T} \omega_{i,t}^d \cdot g_{i,t}^d(x^d) - v^d \cdot h^d(x^d)$, $\omega^d$ is the vector composed of $\omega_{i,t}^d$, and $\omega_{i,t}^d$ and $v^d$ are dual variables. In the equivalent sets of KKT conditions, (5-b) and (6-b) are the gradient equilibrium conditions, and (5-c) and (6-c) enforce complementarity slackness. By replacing LLMs with the equivalent sets, the bi-level model takes the standard form of MPEC as follows:

$$(1\text{-a}) \quad (7\text{-a})$$

subject to:

$$(1\text{-d}), (1\text{-e}), (5\text{-a})\text{-}(5\text{-d}), \forall n, (6\text{-a})\text{-}(6\text{-d}). \quad (7\text{-b})$$

### B. MILP Model Formulation

In the MPEC model, complementarity constraints (5-c) and

(6-c) are non-linear. And constraints (3-b), (3-d)-(3-f), (4-b), (4-d), (5-c), (5-d), (6-c) and (6-d) are in form of the left part of (8), which can be linearized using Fortuny-Amat McCarl linearization method by introducing an auxiliary 0-1 variable $u$ and a sufficiently large positive value $M$, as shown in (8).

$$\begin{cases} a,b \geq 0 \\ a \cdot b = 0 \end{cases} \Rightarrow \begin{cases} 0 \leq a \leq M \cdot u \\ 0 \leq b \leq M \cdot (1-u) \end{cases} \quad (8)$$

After using the above linearization technique, the bi-level model is successfully converted into an MILP problem.

## IV. CASE STUDIES

### A. Case Description

The proposed energy storage sharing strategy is tested on an assumed community of 100 households in California and all the households are served by the same DisCo. For simplicity, we neglect the diversity of customer loads and suppose that half of the households in the community are installed with rooftop photovoltaic (PV) panels. The load profile of these households follow the "duck curve" while the loads of the other households are typical (see Fig. 3). The TOU prices [15] and LMPs [16] are as shown in Fig. 4. Besides, 1) the total capacity of energy storage is 800 kWh; 2) the charging and discharging efficiencies are both set as 0.92; 3) the coefficient $k = 0.25 \text{h}^{-1}$; 4) weight factors $\lambda_1$, $\lambda_2$ and $\lambda_3$ are respectively set as 0.8, 6.69 and 1 according to the capacity benefits of power system [17] and the ratio of the average LMP and the average TOU price; 5) we divide a day evenly into 48 time slots, i.e., $T = 48$, and regard the loads within each time slot as a constant.

Based on the above parameter settings, we carry out the simulation under two LMP profiles, i.e., LMP-1 and LMP-2 in Fig. 4 (b), which are quite different but occur actually in California, respectively. All the problems are solved by software package CPLEX [18] and performed on a laptop with an Intel Core i5 processor and 8GB memory.

### B. Results and Analysis

#### 1) Simulation Results

The numerical simulation results are listed in Tables I-III and the load profiles under LMP-1 and LMP-2 are shown in Figs. 5 and 6, respectively. In Tables II and III, customers 1 and 2 are users with typical load profiles and load profiles with PV panel installed, respectively; scenarios 1, 2 and 3 represent the energy storage are controlled only by DisCo, controlled only by customers, and shared among DisCo and customers, respectively. In Figs. 5 and 6, Load-O is the original load profile; Load-D, Load-C and Load-S are load profiles under scenarios 1, 2 and 3, respectively. Besides, the computation time of the final MILP model is about 0.1s, so the strategy is temporally appropriate for practical operation.

#### 2) Analysis of Results under LMP-1

Under the LMP-1, according to Table I, most of the energy storage are assigned to customers, which can also be verified to some extent by Table II where the results of scenario 2 are more similar to the ones of scenario 3 than the ones of scenario 1. Table II also show that both DisCo and customers can gain the cost saving no matter who completely controls the energy storage, i.e., under either scenario 1 or 2. As the result, the competitive relationship between DisCo and customers for the control of energy storage is weakened. On the other hand, it can be observed in Fig. 5 that the difference among Load-D, Load-C and Load-S is small.

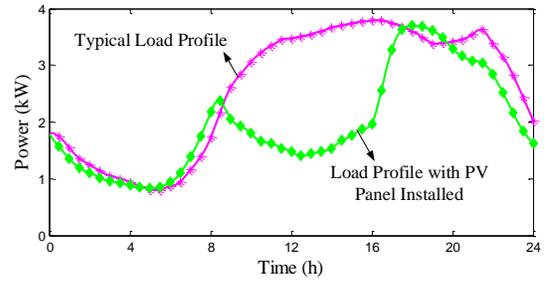

Fig. 3 Load profiles of households.

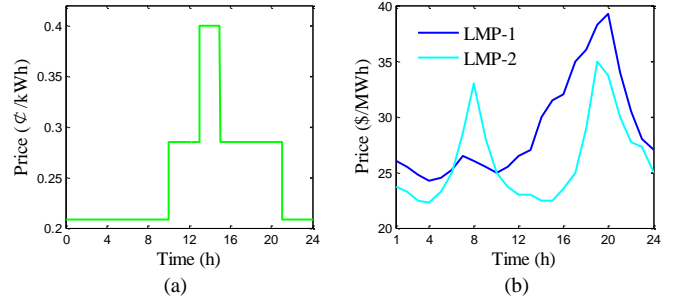

Fig. 4 Curves of (a) TOU prices and (b) LMPs.

TABLE I.  ENERGY STORAGE DIVISION RESULTS (KWH)

| LMP | DisCo | Customers 1 | Customers 2 |
|---|---|---|---|
| LMP-1 | 53 | 261 | 486 |
| LMP-2 | 469 | 116 | 215 |

TABLE II.  ENERGY COST AND PEAK LOAD REDUCTION UNDER LMP-1 (%)

| Scenario | Cost Reduction of DisCo | Cost Reduction of Customers 1 | Cost Reduction of Customers 2 | Peak Load Reduction |
|---|---|---|---|---|
| 1 | 5.28 | 4.21 | 1.84 | 14.29 |
| 2 | 3.82 | 6.24 | 3.61 | 13.89 |
| 3 | 4.08 | 5.72 | 3.48 | 18.26 |

TABLE III.  ENERGY COST AND PEAK LOAD REDUCTION UNDER LMP-2 (%)

| Scenario | Cost Reduction of DisCo | Cost Reduction of Customers 1 | Cost Reduction of Customers 2 | Peak Load Reduction |
|---|---|---|---|---|
| 1 | 5.56 | -4.53 | -2.84 | -3.88 |
| 2 | -2.42 | 6.24 | 3.61 | 13.89 |
| 3 | 3.98 | 4.26 | 2.62 | 13.89 |

Additionally, since the LMP-1 curve and TOU price curve are conforming to the original load profile to some extent, i.e., the price is high when the demand is high, the peak load can be observably shaved under all scenarios 1, 2 and 3 (see Table II and Fig. 5). In spite of this, numerical results demonstrate that the proposed energy storage sharing strategy balances DisCo's and customers' payments and realizes a better performance of peaking shaving (see Table II).

#### 3) Analysis of Results under LMP-2

Different from LMP-1 curve, LMP-2 curve is severely conflicting with TOU price curve and partially conflicting with the original load profile. Thus, Disco's and customers' cost pressure points are not synchronous in time with each

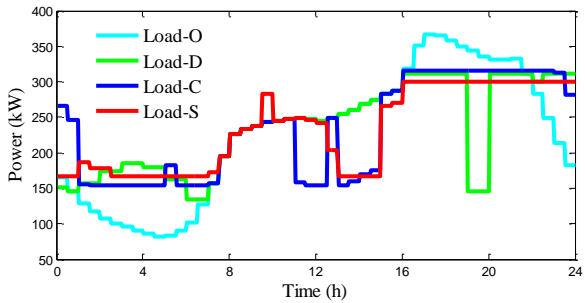
Fig. 5 Load profiles under LMP-1.

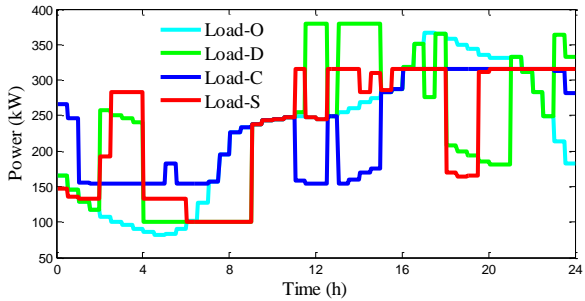
Fig. 6 Load profiles under LMP-2.

other. And due to this, there is no excessive preference in energy storage division between DisCo and customers (see Table I). Also, Load-D and Load-C are in great difference (see Fig. 6).

According to the results in Table III, the conflict of LMP-2 and TOU prices results in higher energy costs of customers and higher peak loads under scenario 1 and higher energy costs of DisCo in scenario 2. However, under scenario 3, the proposed energy storage sharing strategy guarantees that costs of both DisCo and customers are reduced, and effectively restricts the system peak loads.

*4) Summary of Result Analysis*

The proposed energy storage sharing strategy successfully tackle the energy costs trade-off and peak shaving problem among the group of a DisCo and many customers in distribution networks. And the economic loss of DisCo or customers can be well compromised when the conflict of interests exists.

## V. CONCLUSIONS AND FUTURE WORK

In this work, we propose a bi-level model-based energy storage sharing strategy in distribution networks to facilitate the demand response. By using complementarity theory and Fortuny-Amat McCarl linearization method, the model is eventually transformed into an MILP problem without any sacrifice of optimality. And the solving time of the final model is around 0.1s, which provides a basis for practical application of the strategy.

Case studies investigate the performance of the proposed strategy under two kinds of LMP prices. Simulation results demonstrate that the designed sharing of energy storage plays a positive role in the optimization of peak shaving and energy costs of DisCo and customers. Especially when the competitive relationship worsens, the proposed strategy can act as a coordinating part to avoid negative effect on individual energy costs. Hence, our strategy can encourage users and DisCo to participate demand response.

There are some research directions for future work. For instance, we use TOU prices as the electricity prices for customers and do not consider the impacts of different types of electricity retail prices on energy storage sharing, which are worthy of further research. Also, the impacts of various customer load types can be involved in the future work. Another work can be studied is how to charge the DisCo and customers if they join the energy storage sharing plan.